\documentclass[12pt,thmsa]{article}
\usepackage{amsfonts}
\usepackage{amsmath}
\usepackage{amssymb}

\input{tcilatex}
\begin{document}

\title{Abstract Harmonic Analysis on the General Affine Group $GA(n,\mathbb{R%
})$ }
\author{Kahar El-Hussein \\
%EndAName
\textit{Department of Mathematics, Faculty of Science, }\\
\textit{\ Al Furat University, Dear El Zore, Syri\textit{a}\ and}\\
\textit{Department of Mathematics, Faculty of Arts Science Al Quryyat, }\\
\textit{\ Al-Jouf University, KSA }\\
\textit{E-mail : kumath@ju.edu.sa, kumath@hotmail.com }}
\maketitle

\begin{abstract}
Let $GL(n,\mathbb{R})$ be the general linear group and let $GA(n,\mathbb{R})=%
\mathbb{R}^{n}\rtimes _{\rho }GL(n,\mathbb{R})$ be its general affine group.
Let $GL_{+}(n,\mathbb{R})$ be the identity component of $GL(n,\mathbb{R})$,
which consists of the real $n\times n$ matrices with positive determinant
and let $GL_{-}(n,\mathbb{R})$ be the set of all matrices with negative
determinant. Since $GL(n,\mathbb{R})$ is a two copies of $GL_{+}(n,\mathbb{R}%
),$ $i.e$ $GL(n,\mathbb{R})=GL_{+}(n,\mathbb{R})\cup GL_{+}(n,\mathbb{R}),$
because $GL_{-}(n,\mathbb{R})$ has a structure of group isomorphic onto $%
GL_{+}(n,\mathbb{R}),$ see $[13].$ Therefore first we consider the group $%
GA_{+}(n,\mathbb{R})=\mathbb{R}^{n}\rtimes _{\rho }GL_{+}(n,\mathbb{R})$ to
generalize the Fourier transform and to prove the Plancherel theorem for $%
GA_{+}(n,\mathbb{R}).$ Secondly and since $GA(n,\mathbb{R})$ is a two copies
of $GA_{+}(n,\mathbb{R})$, so we can easily to establish the Plancherel
theorem for $GA(n,\mathbb{R}).$
\end{abstract}



\section{\textbf{\ Introduction.}}

\bigskip\ \textbf{1.1. }The general affine group $GA(n,\mathbb{R})=\mathbb{R}%
^{n}\mathbb{\rtimes }_{\rho }$ $GL(n,\mathbb{R})$ of the general linear
group $GL(n,\mathbb{R}),$ which is the semidirect of the real vector $%
\mathbb{R}^{n}$\ with $GL(n,\mathbb{R}).$ Another kind of the affine group
is the Poincare group $\mathbb{R}^{4}\rtimes O(3.1)$, which is the affine
group of the inhomogeneous Lorentz group $O(3.1).$ As well known the affine
groups play an important role in physics in \ cosmology, gauge theory,
gravitation, general relativity, ect... \ .The spacetime symmetry of the
affine modeless of gravity is given by the general affine group $GA(n,%
\mathbb{R})$. The geometrical arena of the theory gravitation and \
electromagnetism is the modified affine frame bundle over a four dimensional
spacetime manifold $M$, the structure group of the frame bundle is the
affine group $GA(4,\mathbb{R})=\mathbb{R}^{4}\mathbb{\rtimes }$ $GL(4,%
\mathbb{R})$. The usual definition of phase space coordinates in terms of
linear frames and use affine frames instead. This leads from $GL(4,\mathbb{R}%
)$ covariance to $GA(4,\mathbb{R})$ covariance and means that the bundle of
linear frames, is replaced by the bundle of affine frames. Also the affine
group $GA(2,\mathbb{R})$ has a fundamental role in the visualization. One
asks can we do the non commutative Fourier analysis on the group $GA(n,%
\mathbb{R})=\mathbb{R}^{n}\mathbb{\rtimes }GL(n,\mathbb{R}).$ In fact that
the abstract harmonic analysis on the locally compact groups is generally a
difficult task. Still now neither the theory of quantum groups nor the
representations theory have done to reach this goal. Recently, these
problems found a satisfactory solution with the papers $[9,10,11].$ In this
paper we will generalize our methods in $[12,13]$ to define the Fourier
transform and to establish the Plancherel theorem for the general affine
group on $GA(n,\mathbb{R}).$

\section{\protect\bigskip Notation and Results for the Nilpotent Lie Group $%
N.$}

\textbf{2.1. }The fine structure of the nilpotent Lie groups will help us to
do the Fourier transform on a connected and simply connected nilpotent Lie
groups $N.$ As well known any group connected and simply connected $N$ has
the following form%
\begin{equation}
N=\left( 
\begin{array}{cccccccccccc}
1 & x_{1}^{1} & x_{1}^{2} & x_{1}^{3} & . & . & . & . & . & x_{1}^{n-2} & 
x_{1}^{n-1} & x_{1}^{n} \\ 
0 & 1 & x_{2}^{2} & x_{2}^{3} & . & . & . & . & . & x_{2}^{n-2} & x_{2}^{n-1}
& x_{2}^{n} \\ 
0 & 0 & 1 & x_{3}^{3} & . & . & . & . & . & x_{3}^{n-2} & x_{3}^{n-1} & 
x_{3}^{n} \\ 
0 & 0 & 0 & 1 & . & . & . & . & . & x_{4}^{n-2} & x_{4}^{n-1} & x_{4}^{n} \\ 
. & . & . & . & . & . & . & . & . & . & . & . \\ 
. & . & . & . & . & . & . & . & . & . & . & . \\ 
. & . & . & . & . & . & . & . & . & . & . & . \\ 
. & . & . & . & . & . & . & . & . & . & . & . \\ 
. & . & . & . & . & . & . & . & . & x_{n-2}^{n-2}. & .x_{n-2}^{n-1} & 
x_{n-2}^{n} \\ 
0 & 0 & 0 & 0 & . & . & . & . & . & 1 & x_{n-1}^{n-2} & x_{n=1}^{n} \\ 
0 & 0 & 0 & 0 & . & . & . & . & . & 0 & 1 & x_{n}^{n} \\ 
0 & 0 & 0 & 0 & . & . & . & . & . & 0 & 0 & 1%
\end{array}%
\right)
\end{equation}

\bigskip As shown, this matrix is formed by the subgroup $\mathbb{R}$, $%
\mathbb{R}^{2}$,...., $\mathbb{R}^{n-1}$, and $\mathbb{R}^{n}$ 
\begin{equation}
\left( \mathbb{R=}\left[ 
\begin{array}{c}
x_{1}^{1} \\ 
1 \\ 
0 \\ 
0 \\ 
. \\ 
. \\ 
. \\ 
. \\ 
. \\ 
0 \\ 
0 \\ 
0%
\end{array}%
\right] ,\mathbb{R}^{2}=\left[ 
\begin{array}{c}
x_{1}^{2} \\ 
x_{2}^{2} \\ 
1 \\ 
0 \\ 
. \\ 
. \\ 
. \\ 
. \\ 
. \\ 
0 \\ 
0 \\ 
0%
\end{array}%
\right] ,..,\mathbb{R}^{n-1}=\left[ 
\begin{array}{c}
x_{1}^{n-1} \\ 
x_{2}^{n-1} \\ 
x_{3}^{n-1} \\ 
x_{4}^{n-1} \\ 
. \\ 
. \\ 
. \\ 
. \\ 
.x_{n-2}^{n-1} \\ 
x_{n-1}^{n-2} \\ 
1 \\ 
0%
\end{array}%
\right] ,\mathbb{R}^{n}=\left[ 
\begin{array}{c}
x_{1}^{n} \\ 
x_{2}^{n} \\ 
x_{3}^{n} \\ 
x_{4}^{n} \\ 
. \\ 
. \\ 
. \\ 
. \\ 
x_{n-2}^{n} \\ 
x_{n=1}^{n} \\ 
x_{n}^{n} \\ 
1%
\end{array}%
\right] \right)
\end{equation}%
$\ \ \ \ \ \ \ \ \ \ \ \ \ \ \ \ \ \ \ \ $

Each $\mathbb{R}^{i}$ is a subgroup of $N$ of dimension $i$ , $1\leq i\leq
n, $ put $d=n+(n-1)+....+2+1=\frac{n(n+1)}{2},$ which is the dimension of $N$
. According to $[9],$ the group $N$ is isomorphic onto the following group

\begin{equation}
(((((\mathbb{R}^{n}\rtimes _{\rho _{n}})\mathbb{R}^{n-1})\rtimes _{\rho
_{n-1}})\mathbb{R}^{n-2}\rtimes _{\rho _{n-2}}.....)\rtimes _{\rho _{2}}%
\mathbb{R}^{2})\rtimes _{\rho _{1}}\mathbb{R}
\end{equation}

\bigskip That means

\begin{equation}
N\simeq (((((\mathbb{R}^{n}\rtimes _{\rho _{n}})\mathbb{R}^{n-1})\rtimes
_{\rho _{n-1}})\mathbb{R}^{n-2}\rtimes _{\rho _{n-2}}.....)\rtimes _{\rho
_{4}}\mathbb{R}^{3}\rtimes _{\rho _{3}}\mathbb{R}^{2})\rtimes _{\rho _{2}}%
\mathbb{R}
\end{equation}

\textbf{2.2}$.$ Denote by $L^{1}(N)$ the Banach algebra of $N$\ that
consists of all complex valued functions on the group $N$, which are
integrable with respect to the Haar measure of $N$ and multiplication is
defined by convolution on $N$ as follows:%
\begin{equation}
g\ast f(X)=\int\limits_{N}f(Y^{-1}X)g(Y)dY
\end{equation}%
for any $f\in L^{1}(N)$ and $g\in L^{1}(N),$ where $X=(X^{1},$ $X^{2},$ $%
X^{3},....,X^{n-2},X^{n-1},X^{n})\in N,$ $%
Y=(Y^{1},Y^{2},Y^{3},....,Y^{n-2},Y^{n-1},Y^{n})\in N,$ $X^{1}=x_{1}^{1},$ $%
X^{2}=(x_{1}^{2},x_{2}^{2}),$ $X^{3}=(x_{1}^{3},x_{2}^{3},x_{3}^{3})$ $%
,...., $ $%
X^{n-2}=(x_{1}^{n-2},x_{2}^{n-2},x_{3}^{n-2},x_{4}^{n-2},...,x_{n-2}^{n-2}),$
$%
X^{n-1}=(x_{1}^{n-1},x_{2}^{n-1},x_{3}^{n-1},x_{4}^{n-1},...,x_{n-2}^{n-1},x_{n-1}^{n-1}), 
$ $%
X^{n}=(x_{1}^{n},x_{2}^{n},x_{3}^{n},x_{4}^{n},...,x_{n-2}^{n},x_{n-1}^{n},x_{n}^{n}) 
$ and $dY=dY^{1}dY^{2}dY^{3},....,dY^{n-2}dY^{n-1}dY^{n}$ is the Haar
measure on $N$ and $\ast $ denotes the convolution product on $N.$ We denote
by $L^{2}(N)$ its Hilbert space. We refer to $[12]$ to define the Fourier
transform on $N$

\textbf{Definition 2.1.} \textit{For} $1\leq i\leq n,$ \textit{let} $%
\mathcal{F}^{i}$\textit{\ be the Fourier transform on} $\mathbb{R}^{i}$ 
\textit{and }$0\leq j\leq n-1,$ \textit{let }$\dprod\limits_{0\leq l\leq j}%
\mathbb{R}^{n-l}=(..((((\mathbb{R}^{n}\rtimes _{\rho _{n}})\mathbb{R}%
^{n-1})\rtimes _{\rho _{n-1}})\mathbb{R}^{n-2}\rtimes _{\rho
_{n-2}})..\times _{\rho _{n-j}}\mathbb{R}^{n-j}),$ \textit{and let }$%
\dprod\limits_{0\leq l\leq j}\mathcal{F}^{n-l}$ $=\mathcal{F}^{n}\mathcal{F}%
^{n-l}\mathcal{F}^{n-2}....\mathcal{F}^{n-j},$ \textit{we can define the
Fourier transform on } $\dprod\limits_{0\leq l\leq n-1}\mathbb{R}^{n-l}=%
\mathbb{R}^{n}\rtimes _{\rho _{n}}\mathbb{R}^{n-1}$ $\rtimes _{\rho _{n-1}}%
\mathbb{R}^{n-2}\rtimes _{\rho _{n-2}}.....\rtimes _{\rho _{3}}\mathbb{R}%
^{2}\rtimes _{\rho _{2}}\mathbb{R}^{1}$\textit{as}%
\begin{eqnarray}
&&\mathcal{F}^{n}\mathcal{F}^{n-1}\mathcal{F}^{n-2}....\mathcal{F}^{2}%
\mathcal{F}^{1}f(\lambda ^{n},\text{ }\lambda ^{n-1},\lambda ^{n-2},\text{%
.......,}\lambda ^{2},\lambda ^{1})  \notag \\
&=&\dint\limits_{N}f(X^{n},X^{n-1},...,X^{2},X^{1})e^{-\text{ }i\langle 
\text{ }(\lambda ^{n},\text{ }\lambda ^{n-1}),(X^{n},X^{n-1}),..,(\lambda
^{2},\text{ }\lambda ^{1}),(X^{2},X^{1})\rangle }\text{ }  \notag \\
&&dX^{n}dX^{n-1}...dX^{2}dX^{1}
\end{eqnarray}%
\textit{for any} $f\in L^{1}(N),$ \textit{where }$%
X=(X^{n},X^{n-1},..,X^{n},X^{n-1}),$ $\mathcal{F}^{d}=\mathcal{F}^{n}%
\mathcal{F}^{n-1}\mathcal{F}^{n-2}....\mathcal{F}^{2}\mathcal{F}^{1}$\ 
\textit{is} \textit{the classical Fourier transform on} $N,$ $%
dX=dX^{n}dX^{n-1}...dX^{2}dX^{1},$ $\lambda =(\lambda ^{n},$ $\lambda
^{n-1},\lambda ^{n-2},..$,$\lambda ^{2},\lambda ^{1}),$ \textit{and}

$\langle (\lambda ^{n},$ $\lambda ^{n-1}),(X^{n},X^{n-1}),..,(\lambda ^{2},$ 
$\lambda ^{1}),(X^{2},X^{1})\rangle =\dsum\limits_{i=1}^{n}X_{i}^{n}\lambda
_{i}^{n}+\dsum\limits_{j=1}^{n-1}X_{j}^{n-1}\lambda
_{j}^{n-1}+..+\dsum\limits_{i=1}^{2}X_{i}^{2}\lambda _{i}^{2}+X^{1}\lambda
^{1}$

\textbf{Plancherels theorem 2.1}. \textit{For any function} $f\in L^{1}(N),$ 
\textit{we have}%
\begin{eqnarray}
&&\dint_{N}\left\vert \mathcal{F}^{d}f(\xi
^{n},X^{n-1},X^{n-2},..,X^{2},X^{1})\right\vert
^{2}dXdX^{n-1}dX^{n-2}..dX^{2}dX^{1}  \notag \\
&&\dint_{N}\left\vert \mathcal{F}^{d}f(\xi ^{n},\lambda ^{n-1},\lambda
^{n-2},..,\lambda ^{2},\lambda ^{1})\right\vert ^{2}d\xi ^{n}d\lambda
^{n-1}d\lambda ^{n-2}..d\lambda ^{2}d\lambda ^{1}
\end{eqnarray}

\bigskip \textit{Proof: }To prove this theorem, we refer to $[12]$.

\section{\protect\bigskip Notation and Results for the Solvable Lie Group $%
AN.$}

\bigskip \textbf{3.1.} Let $G=SL(n,\mathbb{R})$\ be the real semi-simple Lie
group and let $G=KAN$ \ be the Iwasawa decomposition of $G$, where $K=SO(n,%
\mathbb{R}),$and%
\begin{equation}
N=\left( 
\begin{array}{ccccc}
1 & \ast & . & . & \ast \\ 
0 & 1 & \ast & . & \ast \\ 
. & . & . & . & \ast \\ 
. & . & . & . & . \\ 
0 & 0 & . & 0 & 1%
\end{array}%
\right) ,
\end{equation}

\begin{equation}
A=\left( 
\begin{array}{ccccc}
a_{1} & 0 & 0 & . & 0 \\ 
0 & a_{2} & 0 & . & 0 \\ 
. & . & . & . & . \\ 
. & . & . & . & . \\ 
0 & 0 & . & 0 & a_{n}%
\end{array}%
\right)
\end{equation}%
where $a_{1}.a_{2}....a_{n}=1$ and $a_{i}\in \mathbb{R}_{+}^{\ast }.$ The
product $AN$ is a closed subgroup of $G$ and is isomorphic (algebraically
and topologically) to the semi-direct product of $A$ and $N$ with $N$ normal
in $AN.$

Then the group $AN$ is nothing but the group $S=$ $N\rtimes $ $_{\tau }A.$
So the product of two elements $X$ and$Y$ by%
\begin{equation}
(x,\text{ }a)(m,\text{ }b)=(x.\tau (a)y,\text{ }a.b)
\end{equation}%
for\ any $X=(x,a_{1},a_{2},..,a_{n})\in S$ and $Y$ $%
=(m,b_{1},b_{2},..,b_{n})\in S.$ Let $dnda=dmda_{1}da_{2}..da_{n-1}$ be the
right haar measure on $S$ and let $L^{2}(S)$ be the Hilbert space of the
group $S.$ Let $L^{1}(S)$ be the Banach algebra that consists of all complex
valued functions on the group $S$, which are integrable with respect to the
Haar measure of $S$ and multiplication is defined by convolution on $S$ as

\begin{equation}
g\ast f=\int\limits_{S}f((m,b)^{-1}(n,a))g(m,b)dmdb
\end{equation}%
where $dmdb=dmdb_{1}db_{2}..db_{n-1}$ is the right Haar measure on $S=$ $%
N\rtimes $ $_{\rho }A.$ Let $\Lambda =N\times A\times A$ be the group with
law

\begin{equation}
(x,b,a)(y,c,\text{ }d)=(x.\tau (a)y,bc,ad)
\end{equation}%
for any\ $(x,b,a)\in \Lambda ,$ and $(y,c,$ $d)\in \Lambda $

\textbf{Definition 3.1.}\textit{\ For every function }$f$ defined on $S$, 
\textit{one can define a function} $\widetilde{f}$ on $\Lambda $ \textit{as
follows:} 
\begin{equation}
\widetilde{f}(n,a,b)=f(\rho (a)n,ab)
\end{equation}%
\textit{for all} $(n,a,b)\in \Lambda .$ \textit{So every function} $\psi
(n,a)$ \textit{on} $S$\textit{\ extends uniquely as an invariant function} $%
\widetilde{\psi }(n,$ $b,$ $a)$ \textit{on} $\Lambda .$

\ \textbf{Remark 3.1. }\textit{The function} $\widetilde{f}$ \textit{is
invariant in the following sense:} 
\begin{equation}
\widetilde{f}(\rho (s)n,as^{-1},bs)=\widetilde{f}(\rho (s)n,s^{-1}a,bs)=%
\widetilde{f}(n,a,b)
\end{equation}%
\textit{for any} $(n,a,b)\in \Lambda $ \textit{and} $s\in A$

\textbf{Lemma 3.1.}\textit{\ For every function} $f\in L^{1}(S)$ \textit{and
for every }$g\mathcal{\in }$ $L^{1}(S)$, \textit{we have} 
\begin{equation}
g\ast \widetilde{f}(n,a,b)=g\ast _{c}\widetilde{f}(n,a,b)
\end{equation}%
\begin{equation}
\int\limits_{\mathbb{R}^{n-1}}\mathcal{F}_{A}^{n-1}\mathcal{F}^{d}\mathcal{(}%
g\ast \widetilde{f})(\lambda ,\mu ,\nu )d\nu =\mathcal{F}_{A}^{n-1}\mathcal{F%
}^{d}\widetilde{f}(\lambda ,\mu ,0)\mathcal{F}_{A}^{n-1}\mathcal{F}%
^{d}g(\lambda ,\mu )
\end{equation}%
\textit{for every} $(n,a,b)$ $\in \Lambda $, \textit{where} $\ast $\textit{\
signifies the convolution product on} $S$\ \textit{with respect the variables%
} $(n,b)$\ \textit{and} $\ast _{c}$\textit{signifies the commutative
convolution product on} $\ $\textit{the commutative group} $B=N\times A$ 
\textit{with respect the variables }$(n,a),$ \textit{where }$\mathcal{F}%
_{A}^{n-1}$ is the Fourier transform on $A.$

\textbf{Plancherel} \textbf{theorem for }$S.$\textbf{\ 3.1. }\textit{For any
function }$\Psi \in L^{1}(S),$ \textit{we have \ }%
\begin{equation}
\int\limits_{\mathbb{R}^{d}}\int\limits_{\mathbb{R}^{n-1}}\left\vert 
\mathcal{F}^{d}\mathcal{F}_{A}^{n-1}\Psi (\lambda ,\mu )\right\vert
^{2}d\lambda d\mu =\int\limits_{AN}\left\vert \Psi (X,a)\right\vert ^{2}%
\text{ }dXda
\end{equation}

Go back to $[12],$ for the prove of this theorem

\section{\protect\bigskip Fourier Transform and Plancherel Theorem on $GL(n,%
\mathbb{R}).$}

\bigskip \textbf{4.1.} In the following we use the Iwasawa decomposition of $%
G=SL(n,\mathbb{R}),$ to define the Fourier transform and to get Plancherel
theorem on $G=SL(n,\mathbb{R}).$ We denote by $L^{1}(G)$ the Banach algebra
that consists of all complex valued functions on the group $G$, which are
integrable with respect to the Haar measure of $G$ and multiplication is
defined by convolution on $G$ , 
\begin{equation}
\phi \ast f(g)=\int\limits_{G}f(h^{-1}g)\phi (g)dg
\end{equation}

Let $G=SL(n,\mathbb{R})=KNA$ be the Iawsawa decomposition of $G.$The Haar
measure $dg$ on $G$ can be calculated from the Haar measures $dn;$ $da$ and $%
dk$ on $N;A$ and $K;$ respectively, by the formula%
\begin{equation}
\int\limits_{G}f(g)dg=\int\limits_{A}\int\limits_{N}\int%
\limits_{K}f(ank)dadndk
\end{equation}

Keeping in mind that $a^{-2\rho }$ is the modulus of the automorphism $%
n\rightarrow $ $ana^{-1}$ of $N$ we get also the following representation of 
$dg$ 
\begin{equation}
\int\limits_{G}f(g)dg=\int\limits_{A}\int\limits_{N}\int%
\limits_{K}f(ank)dadndk=\int\limits_{N}\int\limits_{A}\int%
\limits_{K}f(nak)a^{-2^{\rho }}dndadk
\end{equation}%
where $\rho =\dim N$ $=\frac{n(n+1)}{2}=1+2+3+....+n-2+n-1+n.$ Furthermore,
using the relation $\int\limits_{G}f(g)dg=\int\limits_{G}f(g^{-1})dg,$ we
receive 
\begin{equation}
\int\limits_{K}\int\limits_{A}\int\limits_{N}f(nak)a^{-2\rho
}dndadk=\int\limits_{K}\int\limits_{A}\int\limits_{N}f(kan)a^{2\rho }dndadk
\end{equation}%
\qquad

\textbf{Definition 4.1.} \textit{The Fourier transform of a function }$f\in
C^{\infty }(K)$\textit{\ is defined as} 
\begin{equation}
Tf(\gamma )=\dint\limits_{K}f(x)\gamma (x^{-1})dx
\end{equation}%
\textit{where }$T$\textit{\ is the Fourier transform on} $K$, \textit{and} \ 
$\gamma \in \widehat{K}$ ($\widehat{K}$ \textit{is the set of irreducible
unitary representations of} $K$ )

\textbf{Theorem (A. Cerezo) 4.1.} \textit{Let} $f\in C^{\infty }(K),$ 
\textit{then we have the inversion of the Fourier transform} 
\begin{equation}
f(x)=\dsum\limits_{\gamma \in \widehat{K}}d\gamma tr[Tf(\gamma )\gamma (x)]
\end{equation}%
\textit{and the Plancherel formula} 
\begin{equation}
\left\Vert f(x)\right\Vert _{2}^{2}=\dint \left\vert f(x)\right\vert
^{2}dx=\dsum\limits_{\gamma \in \widehat{K}}d_{\gamma }\left\Vert Tf(\gamma
)\right\Vert _{H.S}^{2}
\end{equation}%
\textit{for any }$f\in L^{1}(K),$ see $[2,P.562-563],$\textit{where }$%
\left\Vert Tf(\gamma )\right\Vert _{H.S}^{2}$ \textit{is} \textit{the norm
of Hilbert-Schmidt of the operator }$Tf(\gamma ).$

\bigskip \textbf{Definition 4.2}. \textit{For any function} $f\in \mathcal{D}%
(G),$ \textit{we can define a function} $\Upsilon (f)$\textit{on }$G\times K$
\textit{by} 
\begin{equation}
\Upsilon (f)(g,k_{1})=\Upsilon (f)(kna,k_{1})=f(gk_{1})=f(knak_{1})
\end{equation}%
\textit{for }$g=kna\in G,$ \textit{and} $k_{1}\in K$ . \textit{The
restriction of} $\ \Upsilon (f)\ast \psi (g,k_{1})$ \textit{on} $K(G)$ 
\textit{is }$\Upsilon (f)\ast \psi (g,k_{1})\downarrow _{K(G)}=f(nak_{1})\in 
\mathcal{D}(G),$ \textit{and }$\Upsilon (f)(g,k_{1})\downarrow
_{K}=f(g,I_{K})=f(kna)$ $\in \mathcal{D}(G)$

\textbf{Remark 4.1. }\textit{The function }$\Upsilon (f)$ \textit{is
invariant in the following sense }%
\begin{equation}
\Upsilon (f)(gh,h^{-1}k_{1})=f(gk_{1})=f(knak_{1})
\end{equation}

\textbf{Definition 4.3}.\textit{\ Let }$f$ \textit{and }$\psi $ \textit{be
two functions belong to} $\mathcal{D}(G),$ \textit{then we can define the
convolution of } $\Upsilon (f)\ $\textit{and} $\psi $\ \textit{on} $G$ $%
\times K$ \textit{as}

\begin{eqnarray*}
\Upsilon (f)\ast \psi (g,k_{1}) &=&\int\limits_{G}\Upsilon
(f)(gg_{2}^{-1},k_{1})\psi (g_{2})dg_{2} \\
&=&\int\limits_{K}\int\limits_{N}\int\limits_{A}\Upsilon
(f)(knaa_{2}^{-1}n_{2}^{-1}k^{-1}k_{1})\psi
(k_{2}n_{2}a_{2})dk_{2}dn_{2}da_{2}
\end{eqnarray*}

So we get 
\begin{eqnarray*}
\Upsilon (f)\ast \psi (g,k_{1}) &\downarrow &_{K(G)}=\Upsilon (f)\ast \psi
(I_{K}na,k_{1}) \\
&=&\int\limits_{K}\int\limits_{N}\int%
\limits_{A}f(naa_{2}^{-1}n_{2}^{-1}k^{-1}k_{1})\psi
(k_{2}n_{2}a_{2})dk_{2}dn_{2}da_{2}
\end{eqnarray*}%
where $T$ is the Fourier transform on $K,$ and\textit{\ }$I_{K}$ \textit{i}s
the identity element of\textit{\ }$K.$ Denote by\textit{\ }$\mathcal{F}$ is
the Fourier transform on $AN$

\textbf{Definition} \textbf{4.4.} \textit{If }$f\in \mathcal{D}(G)$ \textit{%
and let} $\Upsilon (f)$ \textit{be the associated function to} $f$ , \textit{%
we define the Fourier transform of \ }$\Upsilon (f)(g,k_{1})$ \textit{by }%
\begin{eqnarray}
&&T\mathcal{F}\Upsilon (f))(I_{K},\xi ,\lambda ,\gamma )=T\mathcal{F}%
\Upsilon (f))(I_{K},\xi ,\lambda ,\gamma )  \notag \\
&=&\int_{K}\int_{N}\int_{A}\dsum\limits_{\delta \in \widehat{K}}d_{\delta
}tr[\int_{K}\Upsilon (f)(kna,k_{1})\delta (k^{-1})dk]a^{-i\lambda }e^{-\text{
}i\langle \text{ }\xi ,\text{ }n\rangle }\text{ }\gamma
(k_{1}^{-1})dadndk_{1}  \notag \\
&=&\int_{N}\int_{A}\int_{K}\Upsilon (f)(I_{K}na,k_{1})a^{-i\lambda }e^{-%
\text{ }i\langle \text{ }\xi ,\text{ }n\rangle }\text{ }\gamma
(k_{1}^{-1})dadndk_{1}  \notag \\
&=&\int_{N}\int_{A}\int_{K}f(nak_{1})a^{-i\lambda }e^{-\text{ }i\langle 
\text{ }\xi ,\text{ }n\rangle }\text{ }\gamma (k_{1}^{-1})dadndk_{1}
\end{eqnarray}

\textbf{Theorem 4.2.} (\textbf{\textit{Plancherel's Formula for the Group }}$%
G=SL(n,\mathbb{R}))$\textbf{\textit{. }}\textit{For any function\ }$f\in $%
\textit{\ }$L^{1}(G)\cap $\textit{\ }$L^{2}(G),$\textit{we get }%
\begin{eqnarray}
\int \left\vert f(g)\right\vert ^{2}dg &=&\int_{K}\int_{N}\int_{A}\left\vert
f(kna)\right\vert ^{2}dadndk  \notag \\
&=&\sum_{\gamma \in \widehat{K}}d_{\gamma }\int\limits_{\mathbb{R}%
^{d}}\int\limits_{\mathbb{R}^{n-1}}\left\Vert T\mathcal{F}f(\lambda ,\xi
,\gamma )\right\Vert _{H.S}^{2}d\lambda d\xi  \notag \\
&=&\int\limits_{\mathbb{R}^{d}}\int\limits_{\mathbb{R}^{n-1}}\sum_{\gamma
\in \widehat{K}}d_{\gamma }\left\Vert T\mathcal{F}f(\lambda ,\xi ,\gamma
)\right\Vert _{H.S}^{2}d\lambda d\xi
\end{eqnarray}%
\textit{where} $I_{A},I_{N}$ and $I_{K}$ \textit{are the identity elements of%
} $A$, $N$ \textit{and }$K$ \textit{respectively, }$\mathcal{F}$ \textit{is
the Fourier transform on }$AN$ \textit{and }$T$ \textit{is the Fourier
transform on} $K$

For the proof of this theorem see $[13]$

\textbf{New Group. }Let $GL(n,\mathbb{R)}$ be the general linear group
consisting of all matrices of the form 
\begin{equation}
GL(n,\mathbb{R)}=\{X=\left( 
\begin{array}{c}
a_{ij}%
\end{array}%
\right) ,a_{ij}\in \mathbb{R},\text{ }\prec 1\prec n,,\text{ }\prec j\text{ }%
\prec n,\text{ }and\text{ }\det A\neq 0\}
\end{equation}

As a manifold, $GL(n,\mathbb{R})$ is not connected but rather has two
connected components: the matrices with positive determinant and the ones
with negative determinant which is denoted by $GL_{-}(n,\mathbb{R})$. The
identity component, denoted by $GL_{+}(n,\mathbb{R})$, consists of the real $%
n\times n$ matrices with positive determinant. This is also a Lie group of
dimension $n^{2}$; it has the same Lie algebra as $GL(n,\mathbb{R})$.

The group $GL(n,\mathbb{R})$ is also noncompact. The maximal compact
subgroup of $GL(n,\mathbb{R})$ is the orthogonal group $O(n)$, while the
maximal compact subgroup of $GL_{+}(n,\mathbb{R})$ is the special orthogonal
group $SO(n)$. As for $SO(n)$, the group $GL_{+}(n,\mathbb{R})$ is not
simply connected

\textbf{Theorem 4.3. }$GL_{-}(n,\mathbb{R)}$ \textit{is group isomorphic
onto }$GL_{+}(n,\mathbb{R})$

For the proof of this theorem see $[13]$.

As well known the group $GL_{+}(n,\mathbb{R)}$ is isomorphic onto the direct
product of the two groups $SL(n,\mathbb{R)}$ and $\mathbb{R}_{+}^{\ast },$ $%
i.e$ $GL_{+}(n,\mathbb{R)=}SL(n,\mathbb{R)\times }$ $\mathbb{R}_{+}^{\ast }$
and $GL(n,\mathbb{R)=}GL_{-}(n,\mathbb{R)}\ \cup GL_{+}(n,\mathbb{R)=(}SL(n,%
\mathbb{R)\times }$ $\mathbb{R}_{+}^{\ast })\cup \mathbb{(}SL(n,\mathbb{%
R)\times }$ $\mathbb{R}_{+}^{\ast }).$ Our aim result is

\textbf{Plancherel theorem 4.4. }\textit{Let} $\mathcal{F}_{+}^{\ast }$ 
\textit{be the Fourier transform on }$GL_{+}(n,\mathbb{R}$ \textit{the we get%
} 
\begin{eqnarray}
\int_{GL_{+}(n,\mathbb{R)}}\left\vert f(g,t)\right\vert ^{2}dg\frac{dt}{t}
&=&\int\limits_{\mathbb{R}_{+}^{\ast }}\int_{K}\int_{N}\int_{A}\left\vert
f(kna,t)\right\vert ^{2}dadndk\frac{dt}{t}  \notag \\
&=&\sum_{\gamma \in \widehat{K}}d_{\gamma }\int\limits_{\mathbb{R}%
}\int\limits_{\mathbb{R}^{d}}\int\limits_{\mathbb{R}^{n-1}}\left\Vert 
\mathcal{F}_{+}^{\ast }T\mathcal{F}f(\lambda ,\xi ,\gamma ,\eta )\right\Vert
_{H.S}^{2}d\lambda d\xi d\eta  \notag \\
&=&\int\limits_{\mathbb{R}}\int\limits_{\mathbb{R}^{d}}\int\limits_{\mathbb{R%
}^{n-1}}\sum_{\gamma \in \widehat{K}}d_{\gamma }\left\Vert \mathcal{F}%
_{+}^{\ast }T\mathcal{F}f(\lambda ,\xi ,\gamma ,\eta )\right\Vert
_{H.S}^{2}d\lambda d\xi d\eta
\end{eqnarray}%
for any $f\in L^{1}(GL_{+}(n,\mathbb{R)}\cap L^{2}(GL_{+}(n,\mathbb{R)}$.
For the proof of this theorem see $[13]$.

\textbf{Corollary 4.1. }\textit{Let }$f$ be a function belongs $L^{1}(GL(n,%
\mathbb{R)}\cap L^{2}(GL(n,\mathbb{R)}$ 
\begin{eqnarray}
&&\int_{GL(n,\mathbb{R)}}\left\vert f(g,t)\right\vert ^{2}dg\frac{dt}{t} 
\notag \\
&=&\int_{GL_{-}(n,\mathbb{R)}\cup GL_{+}(n,\mathbb{R)}}\left\vert
f(g,t)\right\vert ^{2}dg\frac{dt}{t}=2\int_{GL_{+}(n,\mathbb{R)}}\left\vert
f(g,t)\right\vert ^{2}dg\frac{dt}{t}  \notag \\
&=&2\int\limits_{\mathbb{R}_{+}^{\ast }}\int_{K}\int_{N}\int_{A}\left\vert
f(kna,t)\right\vert ^{2}dadndk\frac{dt}{t}  \notag \\
&=&2\sum_{\gamma \in \widehat{K}}d_{\gamma }\int\limits_{\mathbb{R}%
}\int\limits_{\mathbb{R}^{d}}\int\limits_{\mathbb{R}^{n-1}}\left\Vert 
\mathcal{F}_{+}^{\ast }T\mathcal{F}f(\lambda ,\xi ,\gamma ,\eta )\right\Vert
_{H.S}^{2}d\lambda d\xi d\eta  \notag \\
&=&2\int\limits_{\mathbb{R}}\int\limits_{\mathbb{R}^{d}}\int\limits_{\mathbb{%
R}^{n-1}}\sum_{\gamma \in \widehat{K}}d_{\gamma }\left\Vert \mathcal{F}%
_{+}^{\ast }T\mathcal{F}f(\lambda ,\xi ,\gamma ,\eta )\right\Vert
_{H.S}^{2}d\lambda d\xi d\eta
\end{eqnarray}

\textbf{Remark 2.1.} this corollary explains the Fourier transform and
Plancherel formula on the non connected Lie group $GL(n,\mathbb{R)}$.

\section{Fourier Transform and Plancherel Theorem for $GA(n,\mathbb{R)}$}

\bigskip \textbf{5.1. }We begin by the group $GA_{+}(n,\mathbb{R})=\mathbb{R}%
^{n}\rtimes _{\rho }GL_{+}(n,\mathbb{R})$ to define the Fourier transform
and to prove the Plancherel formula on the group $GA(n,\mathbb{R}),$where$\
\rho \ $is the group homomorphism from $GL_{+}(n,\mathbb{R})$ into the group 
$Aut(\mathbb{R}^{n})$ of all automorphisms of the real vector group $\mathbb{%
R}^{n}.$ Let $H_{+}=\mathbb{R}^{n}\times GL(n,\mathbb{R})_{+}\times GL_{+}(n,%
\mathbb{R})$\ be the Lie group with the following law 
\begin{equation}
(A,X,Y)(B,P,Q)=(A+\rho (Y)B,XB,YQ)
\end{equation}%
for all $(A,X,Y)\in H_{+}$ , $(C,D,Y)\in H_{+}.$ Denote by $L_{+}=\mathbb{R}%
^{n}\times GL_{+}(n,\mathbb{R})$ the lie group which is direct product of
the two groups $\mathbb{R}^{n}$ and $GL_{+}(n,\mathbb{R}).$ In this case the
group $GA_{+}(n,\mathbb{R})$ can be identified with the closed subgroup $%
\mathbb{R}^{n}\times \left\{ 0\right\} \times _{\rho }GL_{+}(n,\mathbb{R})$
of $H$ $_{+}$ and $L_{+}$ with the subgroup $\mathbb{R}^{n}\times GL_{+}(n,%
\mathbb{R})\times \left\{ 0\right\} $ of $H_{+}$

\ \textbf{Definition 5.1.}\textit{\ For every function }$f$ defined on $%
GA_{+}(n,\mathbb{R})$, \textit{one can define a function} $\widetilde{f}$ on 
$H_{+}$ \textit{as follows:} 
\begin{equation}
\widetilde{f}(A,X,Y)=f(\rho (X)A,XY)
\end{equation}%
\textit{for all} $(A,X,Y)\in H_{+},$ \textit{where }$%
X=(k_{1}n_{1}a_{1},t_{1}),$ $Y=(k_{2}n_{2}a_{2},t_{2}),$ $(k_{1},k_{2})\in
K\times K,$ $(n_{1},n_{2})\in N\times N,$ \textit{and} $(t_{1},t_{2})\in 
\mathbb{R}_{+}^{\ast }\times \mathbb{R}_{+}^{\ast }.$ \textit{So every
function} $\psi (A,Y)\in GA_{+}(n,\mathbb{R})$ \textit{extends uniquely as
an invariant function} $\widetilde{\psi }(A,X,Y)$ \textit{on} $H_{+}.$

\ \textbf{Remark 5.1. }\textit{The function} $\widetilde{f}$ \textit{is
invariant in the following sense:}

\begin{equation}
\widetilde{f}(\rho (T)A,XT^{-1},TY)=\widetilde{f}(A,X,Y)
\end{equation}%
\textit{for any} $(A,X,Y)\in H_{+}$ \textit{and} $T\in GL_{+}(n,\mathbb{R})$

\textbf{Lemma 5.1.}\textit{\ For every function} $f\in L^{1}(GA_{+}(n,%
\mathbb{R}))$ \textit{and for every }$g\mathcal{\in }$ $L^{1}(GA_{+}(n,%
\mathbb{R})$, \textit{we have} 
\begin{equation}
g\ast \widetilde{f}(A,X,Y)=\widetilde{f}\ast _{c}g(A,X,Y)
\end{equation}%
\textit{for every} $(A,X,Y)$ $\in H_{+}$, \textit{where} $\ast $\textit{\
signifies the convolution product on} $GA_{+}(n,\mathbb{R})$\textit{with
respect the variables} $(A,Y)$\ \textit{and} $\ast _{c}$\textit{signifies
the convolution product on} $L_{+}=\mathbb{R}^{n}\times GL_{+}(n,\mathbb{R})$
\textit{with respect the variables }$(A,X),$

\bigskip \textit{Proof: }In fact for any $f$ and $g$ belong $L^{1}(GA_{+}(n,%
\mathbb{R}))$ we have\textit{\ }\ 
\begin{eqnarray}
&&g\ast \widetilde{f}(A,X,Y)=\int\limits_{\mathbb{R}^{n}}\int%
\limits_{GL_{+}(n,\mathbb{R})}\widetilde{f}((B,Q)^{-1}(A,X,Y))g(B,Q)dBdQ 
\notag \\
&=&\int\limits_{\mathbb{R}^{n}}\int\limits_{GL_{+}(n,\mathbb{R})}\widetilde{f%
}\left[ (\rho (Q^{-1})(-B),Q^{-1})(A,X,Y)\right] g(B,Q)dBdQ  \notag \\
&=&\int\limits_{\mathbb{R}^{n}}\int\limits_{GL_{+}(n,\mathbb{R})}\widetilde{f%
}\left[ \rho (Q^{-1})(A-B),X,Q^{-1}Y\right] g(B,Q)dBdQ  \notag \\
&=&\int\limits_{\mathbb{R}^{n}}\int\limits_{GL_{+}(n,\mathbb{R})}\widetilde{f%
}\left[ A-B,XQ^{-1},Y\right] g(B,Q)dBdQ  \notag \\
&=&\widetilde{f}\ast _{c}g(A,X,Y)
\end{eqnarray}

\bigskip \textbf{Definition 5.1. }\textit{If } $f\in L^{1}(GA_{+}(n,\mathbb{R%
})),$ \textit{one can define the Fourier transform of }$f$ \textit{as}%
\begin{eqnarray*}
&&\mathcal{F}_{\mathbb{R}^{n}}\mathcal{F}_{GL_{+}}f(\mu ,(\gamma ,\xi
,\lambda ,\eta )) \\
&=&\int\limits_{\mathbb{R}^{n}}\int_{\mathbb{R}_{+}^{\ast
}}\int_{N}\int_{A}\int_{K}f(A,ank,t)e^{-\text{ }i\langle \text{ }\mu ,\text{ 
}A\rangle }t^{-i\eta }a^{-i\lambda }e^{-\text{ }i\langle \text{ }\xi ,\text{ 
}n\rangle }\text{ }\gamma (k^{-1})dAdadndk\frac{dt}{t}
\end{eqnarray*}%
\textit{here} $\mathcal{F}_{GL_{+}}$ \textit{denotes the Fourier transform on%
} $GL_{+}(n,\mathbb{R})$ \textit{and Let}\textbf{\ }$\mathcal{F}_{\mathbb{R}%
^{n}}$ \textit{be the Fourier transform on} $\mathbb{R}^{n}$

\textbf{Corollary 5.1. }\textit{For any function }$f\in L^{1}(GA_{+}(n,%
\mathbb{R})),$ \textit{we have}%
\begin{eqnarray}
&&\int\limits_{\mathbb{R}}\int\limits_{\mathbb{R}^{d}}\int\limits_{\mathbb{R}%
^{n-1}}\sum_{\delta \in \widehat{K}}d_{\delta }tr[\mathcal{F}_{\mathbb{R}%
^{n}}\mathcal{F}_{GL_{+}}\mathcal{(}g\ast \widetilde{f})(\mu ,(\gamma ,\xi
,\lambda ,,\eta ),(\delta ,\nu ,\sigma ,\omega )]d\nu d\sigma d\omega  \notag
\\
&=&[\mathcal{F}_{\mathbb{R}^{n}}\mathcal{F}_{GL_{+}}\mathcal{(}\widetilde{f}%
)(\mu ,(\gamma ,\xi ,\lambda ,,\eta ),I_{GL_{+}})\mathcal{F}_{\mathbb{R}^{n}}%
\mathcal{F}_{GL_{+}}\mathcal{(}g)(\mu ,(\gamma ,\xi ,\lambda ,\eta ))]
\end{eqnarray}%
\textbf{\ }

\bigskip The proof of this theorem results immediately from lemma \textbf{4.1%
}

\textbf{Plancherel} \textbf{Theorem 5.1. }\textit{For any function }$f\in
L^{1}(GA_{+}(n,\mathbb{R})),$ \textit{we have}%
\begin{eqnarray}
&=&\int\limits_{\mathbb{R}^{n}}\int\limits_{GL_{+}(n,\mathbb{R})}\left\vert
f(B,Q)\right\vert ^{2}dBdQ  \notag \\
&=&\int_{GA_{+}(n,\mathbb{R)}}\left\vert f(A,g,t)\right\vert ^{2}dBdg\frac{dt%
}{t}  \notag \\
&=&\int\limits_{\mathbb{R}^{n}}\int\limits_{\mathbb{R}_{+}^{\ast
}}\int_{K}\int_{N}\int_{A}\left\vert f(A,kna,t)\right\vert ^{2}dBdadndk\frac{%
dt}{t}  \notag \\
&=&\sum_{\gamma \in \widehat{K}}d_{\gamma }\int\limits_{\mathbb{R}%
^{n}}\int\limits_{\mathbb{R}}\int\limits_{\mathbb{R}^{d}}\int\limits_{%
\mathbb{R}^{n-1}}\left\Vert \mathcal{F}_{\mathbb{R}^{n}}\mathcal{F}%
_{GL_{+}}f(\mu ,(\gamma ,\xi ,\lambda ,\eta ))\right\Vert _{H.S}^{2}d\mu
d\lambda d\xi d\eta  \notag \\
&=&\int\limits_{\mathbb{R}^{n}}\int\limits_{\mathbb{R}}\int\limits_{\mathbb{R%
}^{d}}\int\limits_{\mathbb{R}^{n-1}}\sum_{\gamma \in \widehat{K}}d_{\gamma
}\left\Vert \mathcal{F}_{\mathbb{R}^{n}}\mathcal{F}_{GL_{+}}f(\mu ,(\gamma
,\xi ,\lambda ,\eta )\right\Vert _{H.S}^{2}d\mu d\lambda d\xi d\eta
\end{eqnarray}

\textit{Proof}: Let $\widetilde{\overset{\vee }{f}}$ the function defined as 
\begin{equation}
\widetilde{\overset{\vee }{f}}(A,X,Y)=\overset{\vee }{f}(\rho (X)A,XY)=%
\overline{f[(\rho (X)A,XY)^{-1}]}
\end{equation}

Then we have\textit{\ }%
\begin{eqnarray}
&&f\ast \widetilde{\overset{\vee }{f}}(0,I_{GL_{+}},I_{GL_{+}})  \notag \\
&=&\int_{GA_{+}(n,\mathbb{R)}}\widetilde{\overset{\vee }{f}}%
[(B,Q)^{-1}(0,I_{GL_{+}},I_{GL_{+}})]f(B,Q)dBdQ  \notag \\
&=&\int\limits_{\mathbb{R}^{n}}\int\limits_{GL_{+}(n,\mathbb{R})}\widetilde{%
\overset{\vee }{f}}\left[ \rho (Q^{-1})(-B),I_{GL_{+}},Q^{-1}I_{GL_{+}}%
\right] f(B,Q)dBdQ  \notag \\
&=&\int\limits_{\mathbb{R}^{n}}\int\limits_{GL_{+}(n,\mathbb{R})}\widetilde{%
\overset{\vee }{f}}\left[ -B,Q^{-1},I_{GL_{+}}\right] f(B,Q)dBdQ  \notag \\
&=&\int\limits_{\mathbb{R}^{n}}\int\limits_{GL_{+}(n,\mathbb{R})}\overset{%
\vee }{f}\left[ \rho (Q^{-1})(-B),Q^{-1}\right] f(B,Q)dBdQ  \notag \\
&=&\int\limits_{\mathbb{R}^{n}}\int\limits_{GL_{+}(n,\mathbb{R})}\overline{%
f(B,Q)f}(B,Q)dBdQ  \notag \\
&=&\int\limits_{\mathbb{R}^{n}}\int\limits_{GL_{+}(n,\mathbb{R})}\left\vert
f(B,Q)\right\vert ^{2}dBdQ  \notag \\
&=&\int\limits_{\mathbb{R}^{n}}\int_{GA_{+}(n,\mathbb{R)}}\left\vert
f(B,g,t)\right\vert ^{2}dBdg\frac{dt}{t}  \notag \\
&=&\int\limits_{\mathbb{R}^{n}}\int_{GA_{+}(n,\mathbb{R)}}\left\vert
f(B,ank,t)\right\vert ^{2}dBdndadk\frac{dt}{t} \\
&=&\int\limits_{\mathbb{R}^{n}}\int_{GA_{+}(n,\mathbb{R)}}\left\vert
f(B,kna,t)\right\vert ^{2}dBdndadk\frac{dt}{t}
\end{eqnarray}%
where $I_{GL_{+}}$\ is the identity element of $GL_{+}(n,\mathbb{R)}$.\ In
other hand, we get

\begin{eqnarray*}
&&f\ast \widetilde{\overset{\vee }{f}}(0,I_{GL_{+}},I_{GL_{+}}) \\
&=&\int\limits_{\mathbb{R}^{n}}\int\limits_{\mathbb{R}}\int\limits_{\mathbb{R%
}^{d}}\int\limits_{\mathbb{R}^{n-1}}\int\limits_{\mathbb{R}}\int\limits_{%
\mathbb{R}^{d}}\int\limits_{\mathbb{R}^{n-1}}\sum_{\gamma \in \widehat{K}%
}d_{\gamma }\sum_{\delta \in \widehat{K}}d_{\delta }tr[\mathcal{F}_{\mathbb{R%
}^{n}}\mathcal{F}_{GL_{+}}\mathcal{(}f\ast \widetilde{\overset{\vee }{f}}%
)(\mu ,(\gamma ,\xi ,\lambda ,,\eta ),(\delta ,\nu ,\sigma ,\omega )] \\
&&d\mu d\lambda d\xi d\eta d\nu d\sigma d\omega \\
&=&\int\limits_{\mathbb{R}^{n}}\int\limits_{\mathbb{R}}\int\limits_{\mathbb{R%
}^{d}}\int\limits_{\mathbb{R}^{n-1}}\sum_{\gamma \in \widehat{K}}d_{\gamma
}tr[\mathcal{F}_{\mathbb{R}^{n}}\mathcal{F}_{GL_{+}}\mathcal{(}\widetilde{%
\overset{\vee }{f}})(\mu ,(\gamma ,\xi ,\lambda ,,\eta ),I_{GL_{+}})]%
\mathcal{F}_{\mathbb{R}^{n}}\mathcal{F}_{GL_{+}}f(\mu ,\lambda ,\xi ,\gamma
,\eta ) \\
&&d\mu d\lambda d\xi d\eta
\end{eqnarray*}%
\ 

Calculate the formula $\mathcal{F}_{\mathbb{R}^{n}}\mathcal{F}_{GL_{+}}%
\mathcal{(}\widetilde{\overset{\vee }{f}})(\mu ,(\gamma ,\xi ,\lambda ,\eta
),I_{GL_{+}})$

\begin{eqnarray*}
&=&\int\limits_{\mathbb{R}^{n}}\int\limits_{\mathbb{R}}\int\limits_{\mathbb{R%
}^{d}}\int\limits_{\mathbb{R}^{n-1}}\mathcal{F}_{\mathbb{R}^{n}}\mathcal{F}%
_{GL_{+}}\mathcal{(}\widetilde{\overset{\vee }{f}})(\mu ,(\mu ,(\gamma ,\xi
,\lambda ,,\eta ),I_{GL_{+}}) \\
&=&\int\limits_{\mathbb{R}^{n}}\int\limits_{\mathbb{R}}\int\limits_{\mathbb{R%
}^{d}}\int\limits_{\mathbb{R}^{n-1}}\int\limits_{\mathbb{R}^{n}}\int_{%
\mathbb{R}_{+}^{\ast }}\int_{N}\int_{A}\int_{K}\widetilde{\overset{\vee }{f}}%
(A,(kna,t),I_{GL_{+}}) \\
&&e^{-\text{ }i\langle \text{ }\mu ,\text{ }A\rangle }t^{-i\eta
}a^{-i\lambda }e^{-\text{ }i\langle \text{ }\xi ,\text{ }n\rangle }\text{ }%
\gamma (k^{-1})dAdkdnda\frac{dt}{t}d\mu d\lambda d\xi d\eta \\
&=&\int\limits_{\mathbb{R}^{n}}\int\limits_{\mathbb{R}}\int\limits_{\mathbb{R%
}^{d}}\int\limits_{\mathbb{R}^{n-1}}\int\limits_{\mathbb{R}^{n}}\int_{%
\mathbb{R}_{+}^{\ast }}\int_{N}\int_{A}\int_{K}\overset{\vee }{f}(\rho
(ank,t)A,(ank.t)) \\
&&e^{-\text{ }i\langle \text{ }\mu ,\text{ }A\rangle }t^{-i\eta
}a^{-i\lambda }e^{-\text{ }i\langle \text{ }\xi ,\text{ }n\rangle }\text{ }%
\gamma (k^{-1})dkdAdadn\frac{dt}{t}d\mu d\lambda d\xi d\eta \\
&=&\int\limits_{\mathbb{R}^{n}}\int\limits_{\mathbb{R}}\int\limits_{\mathbb{R%
}^{d}}\int\limits_{\mathbb{R}^{n-1}}\int\limits_{\mathbb{R}^{n}}\int_{%
\mathbb{R}_{+}^{\ast }}\int_{N}\int_{A}\int_{K}\overline{f((\rho
(ank,t)A,(ank.t))^{-1})} \\
&&e^{-\text{ }i\langle \text{ }\mu ,\text{ }A\rangle }t^{-i\eta
}a^{-i\lambda }e^{-\text{ }i\langle \text{ }\xi ,\text{ }n\rangle }\text{ }%
\gamma (k^{-1})dkdAdadn\frac{dt}{t}d\mu d\lambda d\xi d\eta \\
&=&\int\limits_{\mathbb{R}^{n}}\int\limits_{\mathbb{R}}\int\limits_{\mathbb{R%
}^{d}}\int\limits_{\mathbb{R}^{n-1}}\int\limits_{\mathbb{R}^{n}}\int_{%
\mathbb{R}_{+}^{\ast }}\int_{N}\int_{A}\int_{K}\overline{%
f((-A,(k^{-1}n^{-1}a^{-1},t^{-1}))} \\
&&e^{-\text{ }i\langle \text{ }\mu ,\text{ }A\rangle }t^{-i\eta
}a^{-i\lambda }e^{-\text{ }i\langle \text{ }\xi ,\text{ }n\rangle }\text{ }%
\gamma (k^{-1})dk]dAdadn\frac{dt}{t}d\mu d\lambda d\xi d\eta \\
&=&\int\limits_{\mathbb{R}^{n}}\int\limits_{\mathbb{R}}\int\limits_{\mathbb{R%
}^{d}}\int\limits_{\mathbb{R}^{n-1}}\int\limits_{\mathbb{R}^{n}}\int_{%
\mathbb{R}_{+}^{\ast }}\int_{N}\int_{A}\int_{K}\overline{f((-A,(kna,t))} \\
&&\overline{e^{-\text{ }i\langle \text{ }\mu ,\text{ }A\rangle }t^{-i\eta
}a^{-i\lambda }e^{-\text{ }i\langle \text{ }\xi ,\text{ }n\rangle }}\text{ }%
\gamma ^{\ast }(k^{-1})dkdAdadn\frac{dt}{t}d\mu d\lambda d\xi d\eta \\
&=&\int\limits_{\mathbb{R}^{n}}\int\limits_{\mathbb{R}}\int\limits_{\mathbb{R%
}^{d}}\int\limits_{\mathbb{R}^{n-1}}\overline{\mathcal{F}_{\mathbb{R}^{n}}%
\mathcal{F}_{GL_{+}}\mathcal{(}f)(\mu ,(\gamma ^{\ast },\xi ,\lambda ,\eta ))%
}d\mu d\lambda d\xi d\eta
\end{eqnarray*}

Finally, we get

\begin{eqnarray}
&&f\ast \widetilde{\overset{\vee }{f}}(0,I_{GL_{+}},I_{GL_{+}})  \notag \\
&=&\int\limits_{\mathbb{R}^{n}}\int\limits_{GL_{+}(n,\mathbb{R})}\left\vert
f(A,Q)\right\vert ^{2}dAdQ  \notag \\
&=&\int\limits_{\mathbb{R}^{n}}\int_{GA_{+}(n,\mathbb{R)}}\left\vert
f(A,g,t)\right\vert ^{2}dAdg\frac{dt}{t}  \notag \\
&=&\int\limits_{\mathbb{R}^{n}}\int_{GA_{+}(n,\mathbb{R)}}\left\vert
f(A,ank,t)\right\vert ^{2}dAdndadk\frac{dt}{t}  \notag \\
&=&\int\limits_{\mathbb{R}^{n}}\int_{GA_{+}(n,\mathbb{R)}}\left\vert
f(A,kna,t)\right\vert ^{2}dAdndadk\frac{dt}{t} \\
&=&\int\limits_{\mathbb{R}^{n}}\int\limits_{\mathbb{R}}\int\limits_{\mathbb{R%
}^{d}}\int\limits_{\mathbb{R}^{n-1}}\sum_{\gamma \in \widehat{K}}d_{\gamma
}tr[\mathcal{F}_{\mathbb{R}^{n}}\mathcal{F}_{GL_{+}}\mathcal{(}\widetilde{%
\overset{\vee }{f}})(\mu ,(\gamma ,\xi ,\lambda ,\eta ),I_{GL_{+}})\mathcal{F%
}_{\mathbb{R}^{n}}\mathcal{F}_{GL_{+}}f(\mu ,\gamma ,\xi ,\lambda ,\eta )] 
\notag \\
&&d\mu d\lambda d\xi d\eta  \notag \\
&=&\int\limits_{\mathbb{R}^{n}}\int\limits_{\mathbb{R}}\int\limits_{\mathbb{R%
}^{d}}\int\limits_{\mathbb{R}^{n-1}}\sum_{\gamma \in \widehat{K}}d_{\delta
}tr[\overline{\mathcal{F}_{\mathbb{R}^{n}}\mathcal{F}_{GL_{+}}\mathcal{(}%
f)(\mu ,(\gamma ^{\ast },\xi ,\lambda ,\eta ))}\mathcal{F}_{\mathbb{R}^{n}}%
\mathcal{F}_{GL_{+}}f(\mu ,(\gamma ,\xi ,\lambda ,,\eta ))]  \notag \\
&=&\int\limits_{\mathbb{R}^{n}}\int\limits_{\mathbb{R}}\int\limits_{\mathbb{R%
}^{d}}\int\limits_{\mathbb{R}^{n-1}}\sum_{\gamma \in \widehat{K}}d_{\delta
}\left\Vert \mathcal{F}_{\mathbb{R}^{n}}\mathcal{F}_{GL_{+}}f(\mu ,(\gamma
,\xi ,\lambda ,,\eta )\right\Vert _{H.S}^{2}d\mu d\lambda d\xi d\eta
\end{eqnarray}

Hence the proof of our theorem. Now we can state our final result

\textbf{Theorem 5.2. }\textit{For any function }$f\in L^{1}(GA(n,\mathbb{R}%
)),$ \textit{we get }%
\begin{eqnarray}
&=&\int\limits_{\mathbb{R}^{n}}\int\limits_{GL(n,\mathbb{R})}\left\vert
f(A,Q)\right\vert ^{2}dAdQ  \notag \\
&=&2\int\limits_{\mathbb{R}^{n}}\int\limits_{GL_{+}(n,\mathbb{R})}\left\vert
f(A,Q)\right\vert ^{2}dAdQ \\
&=&2\int_{GA_{+}(n,\mathbb{R)}}\left\vert f(A,g,t)\right\vert ^{2}dAdg\frac{%
dt}{t}  \notag \\
&=&2\int\limits_{\mathbb{R}^{n}}\int\limits_{\mathbb{R}_{+}^{\ast
}}\int_{K}\int_{N}\int_{A}\left\vert f(A,kna,t)\right\vert ^{2}dAdadndk\frac{%
dt}{t}  \notag \\
&=&2\sum_{\gamma \in \widehat{K}}d_{\gamma }\int\limits_{\mathbb{R}%
^{n}}\int\limits_{\mathbb{R}}\int\limits_{\mathbb{R}^{d}}\int\limits_{%
\mathbb{R}^{n-1}}\left\Vert \mathcal{F}_{\mathbb{R}^{n}}\mathcal{F}%
_{GL_{+}}f(\mu ,(\gamma ,\xi ,\lambda ,,\eta ))\right\Vert _{H.S}^{2}d\mu
d\lambda d\xi d\eta  \notag \\
&=&\int\limits_{\mathbb{R}^{n}}\int\limits_{\mathbb{R}}\int\limits_{\mathbb{R%
}^{d}}\int\limits_{\mathbb{R}^{n-1}}2\sum_{\gamma \in \widehat{K}}d_{\gamma
}\left\Vert \mathcal{F}_{\mathbb{R}^{n}}\mathcal{F}_{GL_{+}}f(\mu ,(\gamma
,\xi ,\lambda ,,\eta )\right\Vert _{H.S}^{2}d\mu d\lambda d\xi d\eta
\end{eqnarray}

\textbf{Conclusion. }I believe that the results of this paper will be a
guideline to study the Fourier analysis on non connected locally compact Lie
groups.

\end{document}